\documentclass{amsart}
\usepackage{amssymb}
\usepackage{amsmath}
\begin{document}
\input epsf.sty

\def\phi{{\varphi}}
\def\d{{\rm d}}
\def\D{{\rm D}}
\def\Q{{\mathbb Q}}
\def\Z{{\mathbb Z}}
\def\C{{\mathbb C}}
\def\N{{\mathbb N}}
\def\P{{\mathbb P}}
\def\R{{\mathbb R}}
\def\deg{{\rm deg\,}}
\def\Re{{\rm Re\,}}
\def\Det{{\rm Det}}
\def\dim{{\rm dim\,}}
\def\Gal{{\rm Gal\,}}
\def\PGL{{\rm PGL\,}}
\def\Rat{{\rm Rat\,}}
\def\Aut{{\rm Aut\,}}
\def\St{{\rm St\,}}

\title{On trees covering chains or stars} 
\author{F. Pakovich}
\address{
Department of Mathematics,
Ben Gurion University of the Negev,
P.O.B. 653, Beer Sheva 84105, Israel}
\email{pakovich@cs.bgu.ac.il }
\date{}
\maketitle
\begin{abstract} In this paper, in the context of the ``Dessins d'enfants'' theory, we give 
a combinatorial criterion for a plane tree
to cover a tree from the classes of ``chains" or ``stars''. Besides, we discuss some applications of this result which are related to the arithmetical theory of torsion on curves. 

\end{abstract}

\def\bp{\begin{proposition}}
\def\ep{\end{proposition}}
\def\bt{\begin{theorem}}
\def\et{\end{theorem}}
\def\be{\begin{equation}}
\def\l{\label}
\def\ee{\end{equation}}
\def\bl{\begin{lemma}}
\def\el{\end{lemma}}
\def\bc{\begin{corollary}}
\def\ec{\end{corollary}}
\def\pr{\noindent{\it Proof. }}
\def\note{\noindent{\bf Note. }}
\def\bd{\begin{definition}}
\def\ed{\end{definition}}
\newtheorem{theorem}{Theorem}[section]
\newtheorem{lemma}{Lemma}[section]
\newtheorem{definition}{Definition}[section]
\newtheorem{corollary}{Corollary}[section]
\newtheorem{proposition}{Proposition}[section]

\section{Introduction} In this paper, in the context of the Grothendieck theory of ``Dessins d'enfants'',
we describe necessary and sufficient combinatorial conditions for an $n$-edged plane tree $\lambda$
to cover a $d$-edged tree from the classes of ``chains" or ``stars''
(see Fig. 1). 
Since for a $d$-edged chain (resp. for a $d$-edged star) the 
corresponding Shabat polynomial is equivalent to the $d$-th Chebyshev polynomial 
$T_d(z)$ (resp. to the polynomial $z^d$)
these conditions correspond to the requirement that, after
an appropriate normalisation, 
the Shabat polynomial $P(z)$ corresponding to $\lambda$ 
admits a 
compositional factorisation of the form $P(z)=T_d(\tilde P(z))$ (resp. of the form $P(z)=(\tilde P(z))^d$). 
Our main result was announced with a sketched proof in the note \cite{p2}. Here we give a detailed proof 
and discuss some applications.
\vskip 0.1cm
\medskip
\epsfxsize=8truecm
\centerline{\epsffile{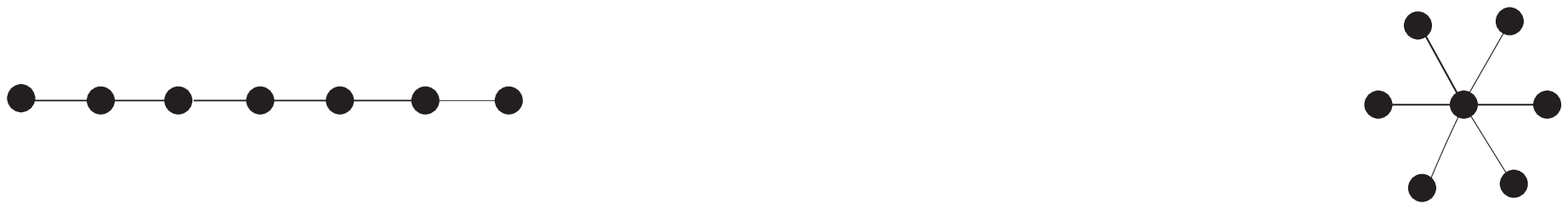}}
\smallskip
\centerline{Figure 1.}
\medskip

For the case of chains the investigated question is related to arithmetics of hyperelliptic curves
via a construction proposed in \cite{p}.  
This construction 
associates to an $n$-edged tree $\lambda$ a hyperelliptic curve $H_{\lambda}$, defined over the field of modules of $\lambda$, such that
the divisor $n(\rho^+_{\infty} -\rho^-_{\infty}),$ where
$\rho^+_{\infty},\rho^-_{\infty}$ are the
points of $H_{\lambda}$ over infinity, is principal.
The order of the divisor $\rho^+_{\infty} -\rho^-_{\infty}$
in the Picard group of $H_{\lambda}$ is equal to $n/d_c,$ where $d_c$ is
a maximal number such that $\lambda$ 
covers a $d_c$-edged chain. This order is an invariant 
with respect to the action of the absolut Galois group $\Gal(\bar\Q/\Q)$ on plane trees and the calculation of this invariant in purely combinatorial terms was the principal motivation for investigations of this paper.

For a tree $\lambda$ define its {\it branch} growing from its vertex
$u$ as a maximal subgraph of $\lambda$ for which 
$u$ is a vertex of valency one. The orientation of the sphere induces
in a natural way a cyclic ordering of branches of $\lambda$ growing from a common vertex. Say that two branches of a tree $\lambda$ 
are {\it adjacent} if they grow from a common
vertex and one of them follows the other with respect to this
ordering. The
number of edges of a branch $a$ is called its {\it weight} and is denoted by $|a|$.

The main result of this paper is the following theorem.

\bt \l{mt}
Let $\lambda$ be an $n$-edged
tree and $d\vert n$. Then $\lambda$ covers a $d$-edged
chain (resp. a $d$-edged star) if and only if the sum
(resp. the difference) of weights of any two
adjacent branches of $\lambda$ is divisible by $d$.
\et

It is not hard to see that for an $n$-edged tree $\lambda$ a number $d_c$ (resp. a number $d_s$) 
such that $\lambda$ 
covers a $d_c$-edged chain (resp. a $d_s$-edged star) is an invariant 
with respect to the action of the group $\Gal(\bar\Q/\Q)$ on trees. Theorem \ref{mt} provides a purely combinatorial description
of these invariants.

\bc \l{co}
For a tree $\lambda$ the invariant $d_c$ (resp. $d_s$) is equal to
the greatest common divisor of all sums $\vert a \vert + \vert b \vert$
(resp. differences $\vert a \vert - \vert b \vert$) such that $a$ and $b$ are adjacent branches of $\lambda.$ 
\ec

The paper has the following structure. At first we
recall a construction from \cite{p} which explains an algebro-geometric meaning of the invariant $d_c$ and discuss some situations
in which Theorem \ref{mt} and Corollary \ref{co} may be useful.
Then we give conditions for a unicellular dessin $\lambda$ to cover an other unicellular dessin, or to be a chain or a star, in terms of arithmetical properties of the canonical involution of oriented edges of $\lambda.$ Finally we prove theorem \ref{mt} and discuss some of its
particular cases.
 
Throughout this paper we will freely use the standard definitions and results of the ``Dessins d'enfants'' theory (see e.g. \cite{sh}, \cite{shl}). Notice that in contrast to the paper \cite{p}
we will assume that all dessins and Belyi functions considered below are clean.

\section{Plane trees and hyperelliptic curves} 

In this subsection we recall a construction from the paper \cite{p} which associates to an $n$-edged tree $\lambda$ with the field of
modules $k_{\lambda}$ a hyperelliptic curve $H_{\lambda}$ defined over $k_{\lambda}$ such that
the divisor $n(\rho^+_{\infty} -\rho^-_{\infty}),$ where
$\rho^+_{\infty},\rho^-_{\infty}$ are
the points of $H_{\lambda}$ over infinity, is principal.

Let $\lambda$ be a tree and let 
$\beta(z)$
be a polynomial from the corresponding equivalence class of Belyi functions.
Set $$H_{\lambda}:
w^2=R(z),$$ where $R(z)$ is a monic polynomial whose (simple) roots are
zeroes of odd multiplicity of the polynomial $\beta(z).$ In other words, if we identify $\lambda$ with the preimage of the segment $[0,1]$ under the map $\beta(z)\,:\, \C\P^1\rightarrow \C\P^1,$ then 
roots of $R(z)$ coincide with vertices of odd valency of $\lambda.$
  
\bp[\cite{p}] \l{prop}
For an $n$-edged tree $\lambda$ the curve
$H_{\lambda}$ is defined over $k_{\lambda}$ and 
the divisor $n(\rho^+_{\infty} -\rho^-_{\infty})$ is principal. 
Furthermore, the order of the divisor $\rho^+_{\infty} -\rho^-_{\infty}$ in the Picard group of $H_{\lambda}$ 
is equal to $n/d_c,$ where $d_c$ is
a maximal number such that $\lambda$ 
covers a $d_c$-edged chain. 
\ep
In order to make Proposition \ref{prop} useful it is important to 
have an expression for the order of the divisor $\rho^+_{\infty} -\rho^-_{\infty}$ in the Picard group of $H_{\lambda}$ in purely combinatorial terms and corollary \ref{co} provides such an expression.
Below we briefly discuss some applications of Proposition \ref{prop}.

For any tree $\lambda$ a total number
of vertices of odd valency  $o_{\lambda}$ is even and $o_{\lambda}=2$ if and only if $\lambda$ is a chain. Furthermore, for 
the genus $g_{\lambda}$ of $H_{\lambda}$ the formula
$g_{\lambda}=(o_{\lambda}-2)/2$ holds. So, the first interesting examples to which the
construction above is applicable are the trees with 4 vertices
of odd valency. This class consists of trees 
homeomorphic either to 
the letter $X$ or to the letter $Y$ (see Fig. 2) and leads to elliptic curves. Note that
after a passage
to the Weierstrass canonical form the divisor $\rho_{\infty}^+-\rho_{\infty}^-$ transforms to a point $(A,B)$ of finite order on  
$H_{\lambda}$ such that $A,B\in k_{\lambda}.$ 

For instance (\cite{p}), the $5$-edged $Y$-tree $\lambda_1$ shown on the Fig. 2 leads to the point $(21, -243)$ of order $5$ on the curve $w^2=4v^3+540v+10665.$ On the other hand, the $6$-edged $X$-tree $\lambda_2$ leads to the point $(3,-16)$ on
the curve $w^2=4v^3+84v-104$. In the last case the order of the corresponding point is 
equal to $3$ since Theorem \ref{mt} implies that $d_c(\lambda_2)=2.$ 
\vskip 0.4cm
\medskip
\epsfxsize=8truecm
\centerline{\epsffile{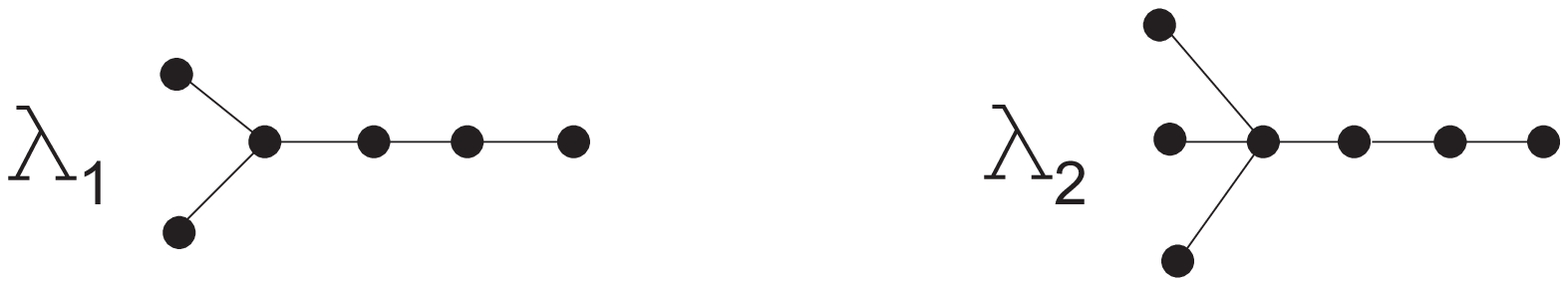}}
\smallskip
\centerline{Figure 2.}
\medskip

Proposition \ref{prop} permits to use for the study of $X$- and
$Y$-trees the well developed arithmetical 
theory of elliptic curves. For example, in the paper \cite{p}, as a corollary of a description of groups $E(\Q)_{tors}$ for elliptic curves over $\Q$ given by 
Mazur \cite{Ma}, {\it a complete list of $Y$-trees defined over $\Q$} was obtained. More general, 
using 
the result of Merel \cite{Me} one can provide a    
{\it lower} bound for the degree of the field of
modules of an $X$- or a $Y$-tree which    
depends only on the invariant $d_c$ (see \cite{p}). 

Another interesting application of 
proposition \ref{prop}
is a method for the finding of
examples of rational divisors of finite order on curves 
defined over $\Q$ (or more generally over number fields) with $g> 1$.
Since for the
curves with $g> 1$ the results similar to the ones of Mazur and Merel do not exist, 
it is interesting how big such an order can be with respect to $g$ (see e.g. \cite{Fl}, \cite{Le}). 
Using proposition \ref{prop} and certain series of
trees
one can obtain for instance
the following result (\cite{p}):
{\it for any
$m$ from the interval $g+1\le m\le 2g+1$ there exists a hyperelliptic
curve of genus $g$ defined over $\Q$ with a rational divisor of order
$m$.} Note that in order to establish this result we do not have to
calculate Belyi functions, all the information needed can be obtained
from the combinatorial analysis of corresponding trees.

\section{Conditions for a unicellular dessin to cover an other
unicellular dessin, or to be a chain or a star} 

Recall that {\it an edge rotation group} $ER(\lambda)$ of a dessin
$\lambda$ is a permutation group of oriented edges of $\lambda$ 
generated by two permutations $\rho_0,$ $\rho_1.$ The permutation  
$\rho_0$ cyclically permutes oriented edges of $\lambda$ around vertices from which they go out
in the order induced by the orientation of the ambient surface, 
and the permutation
$\rho_1$ reverses the orientation of edges. 
Clearly, $ER(\lambda)$ can be identified with the monodromy group of
a Belyi function corresponding to $\lambda.$

Let $\lambda$ be an $n$-edged unicellular dessin. Then the
permutation
$\rho_0\rho_1$ is a cycle of length $2n$.
We associate with $\lambda$ a permutation
$\phi_{\lambda}\subset S_{2n}$ according to the following rule:
enumerate oriented edges of $\lambda$ by the symbols $0,1,\,...\,,2n-1$
in such a way that the cycle $\rho_0\rho_1$ coincides with the cycle
$(01\,...\,2n-1)$ and set $\phi_{\lambda}(i)=
\rho_1(i).$ So, $\phi_{\lambda}$ coincides with
$\rho_1$ but we use a special notation to stress the fact that oriented edges of $\lambda$ are numerated in a specific way. Note that $\phi_{\lambda}$ is a fix points free involution
defined up to a conjugation by some power of the cycle $(01\,...\,2n-1)$.  
Conversely, starting from a fix points free involution $\phi_{\lambda}$ defined on the set $\{0,1,\,...\,,2n-1\}$ we can construct an $n$-edged
unicellular dessin as follows:
enumerate in the counter-clockwise direction the edges of a
$2n$-gon by the numbers $0,1,\,...\,,2n-1$ 
and glue them along $\phi_{\lambda}.$ Two such involutions correspond to the same dessin
if and only if they are conjugated by some power of the cycle $(01\,...\,2n-1)$. 

It is convenient to define the involution 
$\phi_{\lambda}$ on the whole set $\Z$ 
setting the value of $\phi_{\lambda}(j)$, for $j=2nl+\tilde j,$ where $l,\tilde j\in \Z,$ $0\leq \tilde j\leq 2n-1,$
equal to $2nl+\phi_{\lambda}(\tilde j).$

\bp \l{mp}
Let $\lambda$ be an $n$-edged unicellular dessin and 
$d\vert n$. Then $\lambda$ covers a   
$d$-edged dessin $\mu$ if and only if 
\be
\phi_{\lambda}(i+2d)\equiv\phi_{\lambda}(i)\ ({\rm mod}\ 2d)
\ \ \ \ {\it and} \ \ \ \ \phi_{\lambda}(i)\not\equiv i \ ({\rm mod}\ 2d)
\ee for any $i\in \mathbb Z.$
Furthermore, if conditions above are satisfied then $\mu$ is also
unicellular and is defined
uniquely by the condition
$\phi_{\mu}(i)\equiv\phi_{\lambda}(i)\ ({\rm mod}\ 2d)$.
\ep

\pr Indeed, an $n$-edged dessin $\lambda$ covers
a $d$-edged dessin $\mu$ if and only if 
$ER(\lambda)$ has an imprimitivity system $\Omega$ with $2d$ blocks 
such that a permutation induced by $\phi_{\lambda}$ on the set of blocks 
of $\Omega$ has no fix points. Since
$ER(\lambda)$ contains the cycle $(01...2n-1)$ 
such an imprimitivity system should be a collection of the sets
$A_i,$ $0\leq i\leq 2d-1,$ where $A_i$ consists of numbers congruent to $i$ $\rm mod$ $2d.$ Moreover, since the permutations $(01...2n-1)$ and $\rho_1=\phi_{\lambda}$ generate 
$ER(\lambda),$ the collection $A_i,$ $0\leq i\leq 2d-1,$ is an imprimitivity system
if and only if $\phi_{\lambda}(A_i)=A_{\phi(i)}$ for all $i,$ $0\leq i\leq 2d-1.$ 
This condition is equivalent to the first condition of the proposition. The second condition
of the proposition is equivalent to the requirement that 
the permutation of the set of blocks 
of $\Omega$ induced by $\phi_{\lambda}$  
has no fix points.

\bc \l{coro} Let $\lambda$ be an $n$-edged tree. Then $\lambda$ covers a   
$d$-edged tree if and only if 
\be \l{z}
\phi_{\lambda}(i+2d)\equiv\phi_{\lambda}(i)\ ({\rm mod}\ 2d)
\ee for any $i\in \mathbb Z.$
\ec

\pr Indeed, for an $n$-edged tree $\lambda$ we have: 
\be \l{ffb}
\phi_{\lambda}(i)-i \equiv
2\vert a_i \vert -1 \ ({\rm mod}\ 2n),
\ee 
where $a_i$ denotes the branch
of $\lambda$ which contains the oriented edge with number $i$
and grows from the starting point of this edge 
(see Fig. 3).

\medskip
\epsfxsize=13truecm
\centerline{\epsffile{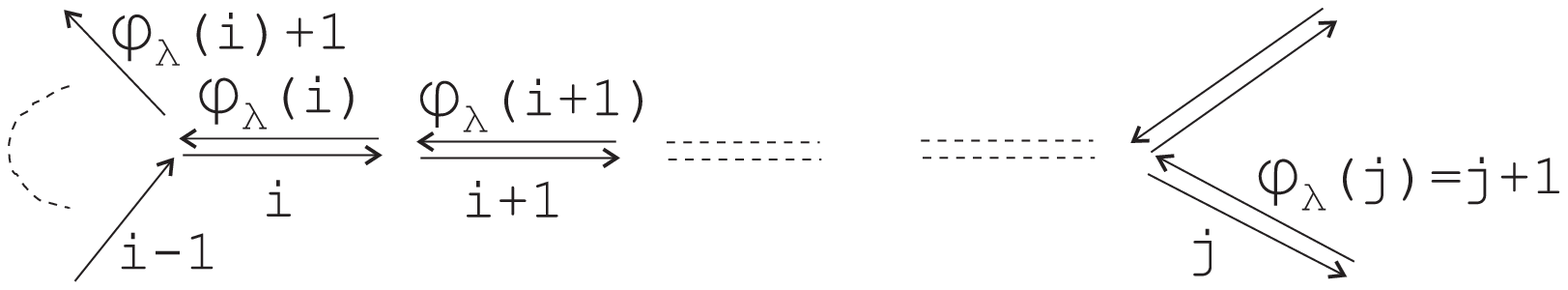}}
\smallskip
\centerline{Figure 3.}

Therefore, the equality \be \l{iop}
\phi_{\lambda}(i)-
i\equiv 1\
({\rm mod}\ 2)
\ee holds and hence the second condition of Proposition \ref{mp} is always satisfied. Furthermore, 
since $\lambda$ is a tree the dessin $\mu$ also is a tree.

\bp \l{pch}
Let $\mu$ be a $d$-edged unicellular dessin. Then $\mu$
is a chain if and only if 
\be \l{ch}
\phi_{\mu}(i)-\phi_{\mu}(i+1)\equiv 1\ ({\rm mod}\ 2d)
\ee for any $i\in \Z.$
\ep

\pr A dessin $\mu$ is a $d$-edged chain if and only if $\phi_
{\mu}
$ 
has the form 
\be \l{c}
\phi_{\mu}(j)=\begin{cases}
\phi_{\mu}(0)-j,&\text{if $0\leq j \leq \phi_{\mu}(0);$}\\
2d+\phi_{\mu}(0)-j,&\text{if $\phi_{\mu}(0) < j \leq 2d-1,$}
\end{cases}
\ee where $\phi_{\mu}(0)$ is an odd number between
$1$ and $2d-1$ (see Fig. 4). Clearly, condition \eqref{c} implies 
condition \eqref{ch}. 

In the opposite direction,
summing equalities \eqref{ch} from $i=0$ to $i=j-1$ we obtain  
$$
\phi_{\mu}(j)\equiv \phi_{\mu}(0)-j \ ({\rm mod}\ 2d).
$$
This implies that $\phi_{\mu}$ has the form \eqref{c}. 
In order to establish that $\phi_{\mu}(0)$ is odd note that
if $\phi_{\mu}(0)=2l$ for some $l,$ $0\leq l \leq d-1,$
then \eqref{c} implies that $\phi_{\mu}(l)=l$ in contradiction with \eqref{iop}.

\bp \l{pst} Let $\mu$ be a $d$-edged unicellular dessin which has at least one vertex of
valency $1.$ Then $\mu$
is a star if and only if 
\be \l{st}
\phi_{\mu}(i)+\phi_{\mu}(i+1)\equiv 2i+1\ ({\rm mod}\
2d)
\ee for any $i\in \Z$. 
\ep

\pr A dessins $\mu$ is a $d$-edged star if and only if $\phi_{\mu}$ 
has the form 
\be \l{s1}
\phi_{\mu}(j)=\begin{cases}
j+(-1)^j\phi_{\mu}(0),&\text{if $0 \leq j+(-1)^j\phi_{\mu}(0) \leq 2d-1;$}\\
j+(-1)^j\phi_{\mu}(0)-2d,&\text{if $2d-1 < j+(-1)^j\phi_{\mu}(0) ;$}\\
2d+j+(-1)^j\phi_{\mu}(0),&\text{if $j+(-1)^j\phi_{\mu}(0) < 0,$}
\end{cases}
\ee 
where either $\phi_{\mu}(0)=1$ or $\phi_{\mu}(0)=2d-1$ (see Fig. 4).

\vskip 0.2cm
\medskip
\epsfxsize=10truecm
\centerline{\epsffile{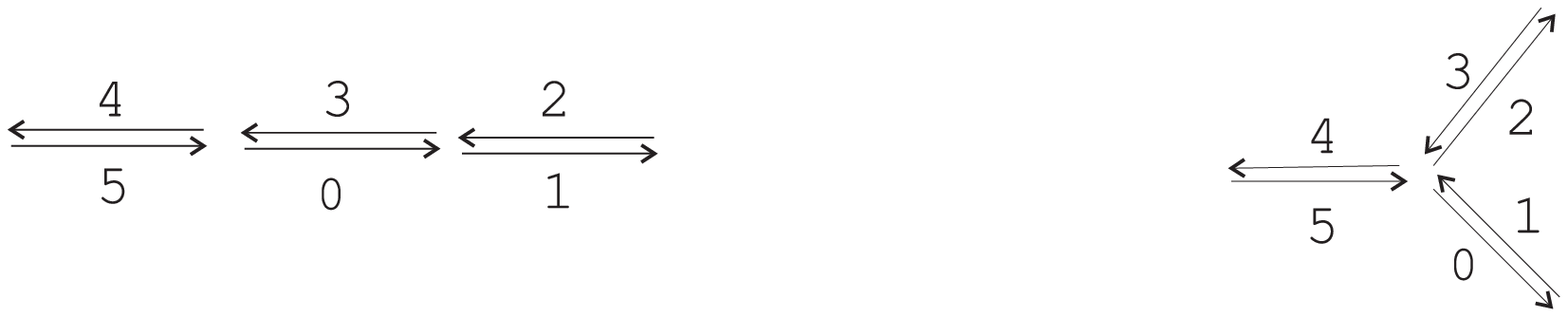}}
\smallskip
\centerline{Figure 4.}
\vskip 0.2cm
\noindent Clearly, condition \eqref{s1} implies condition \eqref{st}. 

In the opposite direction, we have:
$$\phi_{\mu}(0)+(-1)^{j-1}\phi_{\mu}(j)=\sum^{j-1}_{i=0}(-1)^i(\phi_{\mu}(i)+\phi_{\mu}(i+1))\equiv\sum^{j-1}_{i=0}(-1)^i(2i+1) \ ({\rm mod}\
2d).
$$ Since $$\sum^{j-1}_{i=0}(-1)^i(2i+1) 
=(-1)^{j-1} j\,,$$ we conclude that 
\be \l{x}
\phi_{\mu}(j)\equiv j+(-1)^j\phi_{\mu}(0) \ ({\rm mod}\ 2d).
\ee
This implies that $\phi_{\mu}$ has the form \eqref{s1}. 
In order to show that $\phi_{\mu}(0)\equiv\pm 1\ ({\rm mod}\ 2d)
$ note that, since $\mu$ has at least one vertex of 
valency 1, the permutation $\rho_0$ has at least one fixed
point $l$. Since 
$$\rho_0=(01...2d-1)\phi_{\mu},$$ the equalities $\rho_0(l)=l$ and
\eqref{x} imply that
$
\phi_{\mu}(0)\equiv \pm 1 \ ({\rm mod}\
2d).
$
\vskip 0.1cm
\bc \l{lco} An $n$-edged tree covers a   
$d$-edged chain (resp. a   
$d$-edged star) if and only if 
\be \l{ch3}
\phi_{\lambda}(i)-\phi_{\lambda}(i+1)\equiv 1\ ({\rm mod}\ 2d)
\ee 
(resp. 
\be \l{st3}
\phi_{\lambda}(i)+\phi_{\lambda}(i+1)\equiv 2i+1\ ({\rm mod}\
2d)\, )
\ee 
for any $i\in \Z.$
\ec

\pr Indeed, if an $n$-edged tree $\lambda$ covers 
a $d$-edged tree $\mu$ then by corollary \ref{coro} equality \eqref{z} holds. Furthermore, if $\mu$ is a $d$-edged
chain (respectively, a $d$-edged star) then by proposition \ref{pch} (resp. by proposition \ref{pst})  
equality \eqref{ch} (resp. equality \eqref{st}) holds.
Since $\phi_{\mu}(i)\equiv\phi_{\lambda}(i)\ ({\rm mod}\ 2d)$ this implies that condition \eqref{ch3} (resp. condition \eqref{st3}) is satisfied. 

In the opposite direction, arguing as above we conclude that condition \eqref{ch3} (resp. condition \eqref{st3}) implies the condition 
\be \l{rys}
\phi_{\lambda}(j)\equiv \phi_{\lambda}(0)-j \ ({\rm mod}\ 2d)
\ee
(resp. the condition
\be \l{oce}
\phi_{\lambda}(j)\equiv j+(-1)^j\phi_{\lambda}(0) \ ({\rm mod}\ 2d)\, ).
\ee
Since \eqref{rys} as well as \eqref{oce}
implies \eqref{z} it follows from corollary \ref{coro}
that $\lambda$ covers a $d$-edged tree $\mu.$
Furthermore, since \eqref{ch3} (resp. \eqref{st3}) implies  
\eqref{ch} (resp. \eqref{st}), 
proposition \ref{pch} (resp. proposition \ref{pst} taking into account that any tree has vertices of valency one) implies that $\mu$
is a $d$-edged chain (resp. a $d$-edged star).

\section{Proof of theorem \ref{mt}} 
In view of corollary \ref{lco} in order to prove theorem \ref{mt} we only must
show that conditions \eqref{ch3},\eqref{st3} are
actually equivalent to the conditions described in the theorem. 
It follows from formula \eqref{ffb} that for any $i,$ $1\leq i \leq 2n-2,$ we have: 
$$ 
\phi_{\lambda}(i+1)-\phi_{\lambda}(i)+1
=
(\phi_{\lambda}(i+1)-(i+1))+
(i-\phi_{\lambda}(i))+2\equiv $$
$$\equiv 2\vert a_{i+1} \vert +2\vert a_{\phi_{\lambda}(i)} \vert
\ ({\rm mod}\ 2d).
$$ 
Since branches $a_{i+1}$ and $a_{\phi_{\lambda}(i)}$ are adjacent and any adjacent branches have such a form for some $i$ (see Fig. 5), 
this implies that
condition \eqref{ch3} holds if and only if 
the sum 
of the weights of any two
adjacent branches of $\lambda$ is divisible by $d$. 

Similarly,
$$ 
\phi_{\lambda}(i+1)+\phi_{\lambda}(i)-(2i+1)=
(\phi_{\lambda}(i+1)-(i+1))-
(i-\phi_{\lambda}(i))\equiv $$
$$\equiv 2\vert a_{i+1} \vert -2\vert a_{\phi_{\lambda}(i)} \vert
\ ({\rm mod}\ 2d)
$$ and therefore condition \eqref{st3}
is satisfied if and only if the 
difference of the weights of any two
adjacent branches of $\lambda$ is divisible by $d$.

\medskip
\epsfxsize=6truecm
\hskip 0.4cm\centerline{\epsffile{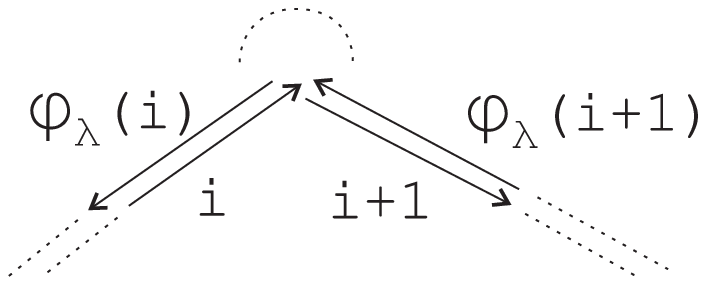}}
\smallskip
\centerline{Figure 5.}
 
\vskip 0.4cm
 
Let us remark 
that the condition of the theorem concerning
chains is automatically satisfied at every vertex of $\lambda$ of valency $2$. In particular, an $n$-edged chain covers a $d$-edged chain
if and only if $d\vert n$ that corresponds to the decomposition $T_n(z)=T_d (T_{n/d}(z))$.
On the
other hand, for a vertex $u$ of odd valency $k$ this condition is equivalent to the requirement that $d\vert \, \vert a \vert$ 
for any branch $a$ growing from $u.$ 
Indeed, the set of all branches growing from $u$ has the form 
$a_{i},a_{\rho_0(i)}, a_{\rho_0^{2}(i)}, ... , a_{\rho_0^{k-1}(i)}$
for some $i,$ $0\leq i \leq 2n-1.$ 
Since $k$ is odd, for any $j\geq 0$ we have:
$$
\vert a_{\rho_0^{j+1}(i)} \vert-\vert a_{\rho_0^{j}(i)}\vert =
\sum_{s=0}^{k-2}(-1)^s (\,\vert a_{\rho_0^{j+s+1}(i)} \vert+\vert a_{\rho_0^{j+s+2}(i)}\vert\ ) \equiv 0 \ ({\rm mod}\ d).$$
Since also 
\be \l{pum}
\vert a_{\rho_0^{j+1}(i)} \vert+\vert a_{\rho_0^{j}(i)}\vert 
\equiv 0 \ ({\rm mod}\
d), \ \ \ \ \ j\geq 0,
\ee 
this implies that $2\vert a_{\rho_0^{j}(i)} \vert \equiv 0 \ ({\rm mod}\
d)$ for any $j,$ $0 \leq  j \leq k-1.$ Hence, 
either $\vert a_{\rho_0^{j}(i)} \vert \equiv 0 \ ({\rm mod}\
d)$ 
or $\vert a_{\rho_0^{j}(i)} \vert \equiv d/2 \ ({\rm mod}\
d).$ 
If for all $j,$ $0 \leq  j \leq k-1,$ we have $\vert a_{\rho_0^{j}(i)} \vert \equiv d/2 \ ({\rm mod}\ d)$  
then summing these equalities and taking into account that $k$ is odd we conclude that
$$n=\vert a_{i} \vert +\vert a_{\rho_0^{}(i)} \vert + ... +\vert  a_{\rho_0^{k-1}(i)} \vert\equiv d/2 \ ({\rm mod}\
d)$$ in contradiction with $d\vert n.$ Therefore, $d \vert\, \vert a_{\rho_0^{j}(i)} \vert $ for at least one $j,$ $0 \leq  j \leq k-1.$ It follows now 
from equalities \eqref{pum} by induction
that $d \vert\, \vert a_{\rho_0^{j}(i)} \vert $ for all $j,$ $0 \leq  j \leq k-1.$

\bibliographystyle{amsplain}

\end{document}